\documentclass[11pt, reqno]{amsart}

\usepackage{amsmath, amsfonts, amssymb, amsbsy, bigstrut, graphicx, enumerate,  upref, longtable, comment}

\usepackage{pstricks}

\newtheorem{thm}{Theorem}[section]

\theoremstyle{remark}

\newcommand{\cov}{\mathrm{Cov}}

\newcommand{\ee}{\mathbb{E}}

\newcommand{\mx}{\mathcal{X}}

\newcommand{\pp}{\mathbb{P}}

\newcommand{\ra}{\rightarrow}
\newcommand{\rr}{\mathbb{R}}

\newcommand{\var}{\mathrm{Var}}

\newcommand{\xp}{X^\prime}

\newcommand{\xx}{\mathcal{X}}

\newcommand{\fpar}[2]{\frac{\partial #1}{\partial #2}}

\newcommand{\mpar}[3]{\frac{\partial^2 #1}{\partial #2 \partial #3}}

\numberwithin{equation}{section}

\usepackage{hyperref}

\begin{document}
\title{A central limit theorem for a new statistic on permutations}
\author{Sourav Chatterjee}
\author{Persi Diaconis}
\address{\newline Departments of Mathematics and Statistics\newline Stanford University \newline \textup{\tt souravc@stanford.edu} \newline
\textup{\tt diaconis@math.stanford.edu}}
\thanks{Sourav Chatterjee's research was partially supported by NSF grant DMS-1441513.}
\thanks{Persi Diaconis's research was partially supported by NSF grant DMS-1208775.}

\begin{abstract}
This paper does three things: It proves a central limit theorem for novel permutation statistics (for example, the number of descents plus the number of descents in the inverse). It provides a clear illustration of a new approach to proving central limit theorems more generally. It gives us an opportunity to acknowledge the work of our teacher and friend B.~V.~Rao.
\end{abstract}

\maketitle

\section{Introduction}
Let $S_n$ be the group of all $n!$ permutations of $\{1,\ldots, n\}$. A variety of statistics $T(\pi)$ are used to enable tasks such as tests of randomness of a time series, comparison of voter profiles when candidates are ranked, non-parametric statistical tests and evaluation of search engine rankings. A basic feature of a permutation is a local `up-down' pattern. Let the number of descents be defined~as 
\[
D(\pi) := |\{i: 1\le i\le n-1,\, \pi(i+1)<\pi(i)\}|\,.
\]
For example, when $n=10$, the permutation $\pi = (\underline{7} \; 1\; \underline{5}\; 3\; \underline{10}\; \underline{8}\; \underline{6}\; 2\; 4\; 9)$ has $D(\pi)= 5$. The enumerative theory of permutations by descents has been intensively studied since Euler. An overview is in Section \ref{overview} below. In seeking to make a metric on permutations using descents we were led to study
\begin{equation}\label{tpi}
T(\pi) := D(\pi)+D(\pi^{-1})\,.
\end{equation}
For a statistician or a probabilist it is natural to ask ``Pick $\pi\in S_n$ uniformly; what is the distribution of $T(\pi)$?'' While a host of limit theorems are known for $D(\pi)$, we found $T(\pi)$ challenging. A main result of this paper establishes a central limit theorem. 
\begin{thm}\label{clt}
For $\pi$ chosen uniformly from the symmetric group $S_n$, and $T(\pi)$ defined by \eqref{tpi}, for $n\ge 2$
\[
\ee(T(\pi)) = n-1\,,\; \; \var(T(\pi)) = \frac{n+7}{6}-\frac{1}{n}\,,
\]
and, normalized by its mean and variance, $T(\pi)$ has a limiting standard normal distribution. 
\end{thm}
The proof of Theorem \ref{clt} uses a method of proving central limit theorems for complicated functions of independent random variables due to Chatterjee~\cite{cha08}. This seems to be a useful extension of Stein's approach. Indeed, we were unable to prove Theorem~\ref{clt} by standard variations of Stein's method such as dependency graphs, exchangeable pairs or size-biased couplings. Theorem \ref{clt} is a special case of the following more general result. Call a statistic $F$ on $S_n$ ``local of degree $k$'' if $F$ can be expressed as
\[
F(\pi) = \sum_{i=0}^{n-k} f_i(\pi)\,,
\] 
where the quantity $f_i(\pi)$ depends only on the relative ordering of $\pi(i+1),\ldots, \pi(i+k)$. For example, the number of descents is local of degree $2$, and the number of peaks is local of degree $3$. We will refer to $f_0,\ldots, f_{n-k}$ as the ``local components'' of $F$. 
\begin{thm}\label{genclt}
Suppose that $F$ and $G$ are local statistics of degree $k$ on $S_n$, as defined above. Suppose further that the absolute values of the local components of $F$ and $G$ are uniformly bounded by $1$. Let $\pi$ be be chosen uniformly from $S_n$ and let 
\[
W := F(\pi) + G(\pi^{-1})\,.
\]
Let $s^2 := \var(W)$. Then the Wasserstein distance between $(W-\ee(W))/s$ and the standard normal distribution is bounded by $C(n^{1/2}s^{-2} + n s^{-3}) k^3$, where $C$ is a universal constant. 
\end{thm}
After the first draft of this paper was posted on arXiv, it was brought to our notice that the joint normality of  $D(\pi)$, $D(\pi^{-1})$ was proved in Vatutin~\cite{vatutin96} in 1996 via a technical tour de force with generating functions. The asymptotic normality in Theorem \ref{clt} follows as a corollary of Vatutin's theorem. Theorem \ref{genclt} is a new contribution of this paper. 

In outline, Section \ref{metrics} describes metrics on permutations and our motivation for the study of $T(\pi)$. Section \ref{overview} reviews the probability and combinatorics of $D(\pi)$ and $T(\pi)$. Section \ref{method} describes Chatterjee's central limit theorem. Section \ref{proof} proves Theorems \ref{clt} and the proof of Theorem \ref{genclt} is in Section \ref{proof2}. Section \ref{future} outlines some other problems where the present approach should work.
\vskip.2in
\noindent{\bf Acknowledgments.} This work derives from conversations with Ron Graham. Walter Stromquist provided the neat formula for the variance of $T(\pi)$. Jason Fulman provided some useful references. Vladimir Vatutin brought the important reference \cite{vatutin96} to our notice. Finally, B.~V.~Rao has inspired both of us by the clarity, elegance and joy that he brings to mathematics.

\section{Metrics on permutations}\label{metrics}
A variety of metrics are in widespread use in statistics, machine learning, probability, computer science and the social sciences. They are used in conjunction with statistical tests, evaluation of election results (when voters rank order a list of candidates), and for combing results of search engine rankings. The book by Marden \cite{marden} gives a comprehensive account of various approaches to statistics on permutations. The book by Critchlow \cite{critchlow} extends the use of metrics on permutations to partially ranked data (top $k$ out of $n$) and other quotient spaces of the symmetric group. The Ph.D.~thesis of Eric Sibony (available on the web) has a comprehensive review of machine learning methods for studying partially ranked data. Finally, the book by Diaconis \cite[Chapter 6]{diaconis} contains many metrics on groups and an extensive list of applications. 

Some widely used metrics are:
\begin{itemize}
\item Spearman's footrule: $d_s(\pi, \sigma) = \sum_{i=1}^n |\pi(i)-\sigma(i)|$.
\item Spearman's rho: $d_\rho^2(\pi, \sigma) = \sum_{i=1}^n (\pi(i)-\sigma(i))^2$. 
\item Kendall's tau: $d_\tau(\pi,\sigma) = $ minimum number of adjacent transpositions to bring $\sigma$ to $\pi$. 
\item Cayley: $d_C(\pi, \sigma) = $ minimum number of transpositions to bring $\sigma$ to $\pi$. 
\item Hamming: $d_H(\pi, \sigma) = |\{i: \pi(i)\ne \sigma(i)\}|$. 
\item Ulam: $d_U(\pi,\sigma) = n - \text{length of longest increasing subsequence in } \pi\sigma^{-1}$.
\end{itemize}
All of these have known means, variances and limit laws \cite{diaconis}. Some of this is quite deep mathematically. For example, the limit law for Ulam's metric is the Tracy--Widom distribution of random matrix theory.

In addition to the metric properties, metrics can be classified by their invariance properties; for example, right invariance ($d(\pi, \sigma) = d(\pi\eta, \sigma\eta)$), left invariance ($d(\pi, \sigma) = d(\eta \pi, \eta \sigma)$), two-sided invariance ($d(\pi, \sigma) = d(\eta_1\pi\eta_2, \eta_1\sigma\eta_2)$), and conjugacy invariance ($d(\pi, \sigma) = d(\eta^{-1}\pi\eta, \eta^{-1}\sigma\eta)$). Common sense requires right invariance; if $\pi$ and $\sigma$ are the rankings of a class on the midterm and on the final, we would not want $d(\pi,\sigma)$ to depend on the class being listed by last name or identity number. All of the metrics above are right invariant. Only the Cayley and Hamming metrics are bi-invariant.
\vskip.1in
\noindent{\underline{A metric from descent structure?}}

It is natural to try and make a metric from descents. Let us call this $d_D(\pi, \sigma)$. By right invariance only the distance $d_D(\text{id}, \sigma)$ must be defined (and then $d_D(\pi, \sigma) = d_D(\text{id}, \sigma\pi^{-1})$). A zeroth try is $d^0_D(\text{id}, \sigma) = D(\sigma)$; at least $d_D^0(\text{id}, \text{id}) = 0$. However, symmetry requires $d_D^0(\text{id}, \sigma) = d^0_D(\sigma, \text{id}) = d_D^0(\text{id}, \sigma^{-1})$ and $D(\sigma)\ne D(\sigma^{-1})$ for many $\sigma$ (for example, when $n=4$, $\sigma = (2\; 4\; 1\; 3)$, $\sigma^{-1}= (3\; 1\; 4\; 2)$, $D(\sigma) = 1$, $D(\sigma^{-1})=2$). In small samples $D(\sigma)=D(\sigma^{-1})$ occurs fairly often. However, in Section \ref{proof} we will prove a bivariate central limit theorem for $D(\sigma)$, $D(\sigma^{-1})$, which suggests that the chance of $D(\sigma)=D(\sigma^{-1})$ is asymptotic to $Cn^{-1/2}$ when $n$ is large. 

A next try is $d_D^1(\text{id}, \sigma) = D(\sigma)+D(\sigma^{-1})$. Then $d_D^1(\pi, \sigma) = D(\sigma \pi^{-1}) + D(\pi\sigma^{-1}) = d_D^1(\sigma, \pi)$. Alas, Ron Graham showed us simple examples where this definition fails to satisfy the triangle inequality! Take $\pi = (3\;4\;1\;2\;5)$, $\sigma =(1\;4\;5\;2\;3)$. A simple check shows that
\[
2+2 = d_D^1(\pi, \text{id}) + d_D^1(\text{id}, \sigma) < d_D^1(\pi, \sigma) = 6\,.
\]
The next idea does work. Form a graph with vertices the $n!$ permutations and an edge from $\pi$ to $\sigma$ with weight $D(\pi\sigma^{-1}) + D(\sigma \pi^{-1})$.  Define $d_D^2(\pi, \sigma)$ as the minimum sum of the weights of paths from $\pi$ to $\sigma$. Experiments show that {\it usually} the minimum path is the edge from $\pi$ to $\sigma$. But the example above shows this is not always the case. We believe that for almost all pairs the graph distance equals the edge weight.

The statistic $T(\pi) = D(\pi)+D(\pi^{-1})$ arose from these considerations. Of course, $T$ does not have to give rise to a metric to be a useful measure of disarray. The Kullback--Leibler `divergence' is a case in point.

\section{Combinatorics and probability for descents}\label{overview}
The study of descents starts with Euler. In studying power series which allow closed form evaluation, Euler showed that 
\begin{equation}\label{euler}
\sum_{k=0}^\infty k^n t^k = \frac{A_n(t)}{(1-t)^{n+1}}\,,
\end{equation}
where the Eulerian polynomial is
\[
A_n(t) = \sum_{\pi\in S_n} t^{D(\pi)+1} = \sum_{i=1}^n A_{n,i} t^i\,,
\]
with $A_{n,i} = |\{\pi\in S_n: D(\pi)=i-1\}|$, the Eulerian numbers. Thus,
\[
\sum_{k=0}^\infty t^k = \frac{1}{1-t}\,,\;\;\sum_{k=0}^\infty kt^k = \frac{t}{(1-t)^2}\,, \;\;\sum_{k=0}^\infty k^2 t^k = \frac{t+t^2}{(1-t)^3}\, , \;\;\ldots
\]
Fulman \cite{fulman99} connects \eqref{euler} and other descent identities to `higher math'.

The Eulerian numbers and polynomials satisfy a host of identities and have neat generating functions. We recommend Carlitz~\cite{carlitz}, Petersen~\cite{petersen13},  Stanley~\cite[p.\ 6]{stanley77}, Graham--Knuth--Patashnik~\cite[Chapter 6]{grahametal}, the book by Petersen~\cite{petersen15} and the references in Sloane~\cite[Seq A008292]{sloane} for basics with pointers to a large literature. 

The probability theory of descents is similarly well studied. The central limit theorem reads:
\begin{thm}\label{descent}
For $\pi$ chosen uniformly in $S_n$, 
\[
\ee(D(\pi)) = \frac{n-1}{2}\,,\;\;\var(D(\pi)) = \frac{n+1}{12}\,,
\]
and, normalized by its mean and variance, $D(\pi)$ has a limiting standard normal distribution.
\end{thm}
\noindent{\it \underline{Remark.}} We point here to six different approaches to the proof of Theorem \ref{descent}. Each comes with an error term (for the Kolmogorov distance) of order $n^{-1/2}$. The first proof uses $m$-dependence: For $1\le i\le n-1$, let 
\[
X_i(\pi) = 
\begin{cases}
1&\text{ if } \pi(i+1)<\pi(i),\\
0 &\text{ else.}
\end{cases}
\]
It is easy to see that $X_1,\ldots, X_{n-1}$ are $2$-dependent. The central limit theorem follows. See Chen and Shao \cite{chenshao} for a version with error terms. A second proof uses a geometrical interpretation due to Stanley~\cite{stanley77}. Let $U_1,\ldots, U_n$ be independent uniform random variables on $[0,1]$. Stanley shows that for all $n$ and $0\le j\le n-1$, 
\[
\pp(D(\pi)=j) = \pp(j < U_1+\cdots +U_n < j+1)\,.
\]
From here the classical central limit theorem for sums of i.i.d.~random variables gives the claimed result. A third proof due to Harper~\cite{harper} uses the surprising fact that the generating function $A_n(t)$ has all real zeros. Now, general theorems for such generating functions show that $D(\pi)$ has the same distribution as the sum of $n$ independent Bernoulli random variables with success probabilities determined by the zeros. Pitman \cite{pitman} surveys this topic and shows
\[
\sup_{-\infty<x<\infty}\left|\pp\left(\frac{D(\pi)-\frac{n-1}{2}}{\sqrt{\frac{n+1}{12}}} \le x\right) - \Phi(x)\right|\le \sqrt{\frac{12}{n}}\,,
\]
where $\Phi$ is the standard normal cumulative distribution function.

A fourth approach due to Fulman \cite{fulman04} uses a clever version of Stein's method of exchangeable pairs (see also Conger \cite{conger}). Each of the papers cited above gives pointers to yet other proofs. A related probabilistic development is in Borodin, Fulman and Diaconis~\cite{borodinetal}. They show that the descent process is a determinantal point process and hence the many general theorems for determinantal point processes apply.

The fifth approach is due to Bender~\cite{bender}, who uses generating functions to prove the CLT. Lastly, David and Barton~\cite{davidbarton} give a proof using the method of moments.

We make two points regarding the above paragraphs. First, all of the proofs depend on some sort of combinatorial magic trick. Second, we were unable to get any of these techniques to work for~$T(\pi)$. 

There have been some applications of descents and related local structures (for example, peaks) in statistics. This is nicely surveyed in Warren and Seneta~\cite{warrenseneta} and Stigler \cite{stigler}.
\vskip.1in
\noindent{\underline{Joint distribution of $D(\pi)$, $D(\pi^{-1})$}}

There have been a number of papers that study the joint distribution of $D(\pi)$, $D(\pi^{-1})$ (indeed, along with other statistics). We mention Rawlings~\cite{rawlings} and Garsia and Gessel~\cite{garsiagessel}. As shown below, these papers derive generating functions in a sufficiently arcane form that we have not seen any way of deriving the information we need from them. 

Two other papers seem more useful. The first by Kyle Petersen \cite{petersen13} treats only $D(\pi)$, $D(\pi^{-1})$ and is very accessible. The second, by Carlitz, Rosellet and Scoville~\cite{carlitzetal}, gives useful recurrences via `manipulatorics'. 

First, let
\[
A_{n,r,s} = |\{\pi\in S_n: D(\pi)=r-1,\, D(\pi^{-1}) = s-1\}| \,.
\]
A table of $A_{8,r,s}$ from Peterson \cite[p.~12]{petersen13} shows a roughly elliptic shape and suggests that the limiting distribution of $D(\pi)$, $D(\pi^{-1})$ might be a bivariate normal distribution with vanishing correlations. Theorem \ref{clt2} from Section \ref{proof} proves normality with limiting correlation zero. 

Carlitz et al.~\cite{carlitzetal} delineate which $r$, $s$ are possible:
\[
A_{n,r,s}=0 \iff r\ge \frac{s-1}{s} n + 1\,.
\]
Let 
\[
A_n(u,v) =  \sum_{\pi\in S_n} u^{D(\pi)+1} v^{D(\pi^{-1})+1} = \sum_{i,j=1}^{n} A_{n,i,j} u^i v^j\,.
\]
Petersen \cite[Theorem 2]{petersen13} gives the following formula for the generating function defined above:
\[
\frac{A_n(u,v)}{(1-u)^{n+1}(1-v)^{n+1}} = \sum_{k,l\ge 0} {kl+n-1\choose n} u^k v^l\,.
\]
From this, he derives a recurrence which may be useful for deriving moments:
\begin{align*}
nA_n(u,v) &= (n^2uv + (n-1)(1-u)(1-v))A_{n-1}(u,v)\\
&\qquad + nuv(1-u)\fpar{}{u}A_{n-1}(u,v) + nuv(1-v) \fpar{}{v} A_{n-1}(u,v)\\
&\qquad + uv(1-u)(1-v) \mpar{}{u}{v} A_{n-1}(u,v)\,.
\end{align*}
From the above identity, he derives a complicated recurrence for $A_{n,i,j}$.  We are interested in $D(\pi)+D(\pi^{-1})$, so the relevant generating function is
\[
A_n(u,u) = (1-u)^{2n+2} \sum_{k,l\ge 0} {kl+n-1\choose n} u^{k+l}\,.
\] 
Finally, Carlitz et al.~\cite{carlitzetal} give
\[
\sum_{n=0}^\infty\sum_{i=0}^\infty \sum_{j=0}^\infty A_{n,i,j} z^ix^jy^n(1-x)^{-(n+1)} (1-z)^{-(n+1)} = \sum_{k=0}^\infty \frac{z^k}{1-x(1-y)^{-k}}\,.
\]
Complicated and seemingly intractable as they are, the generating functions displayed above are in fact amenable to analysis. In a remarkable piece of work from twenty years ago, Vatutin \cite{vatutin96} was able to use these generating functions to prove a class of multivariate central limit theorems for functions of $\pi$ and $\pi^{-1}$. The results were generalized to other settings in Vatutin~\cite{vatutin94} and Vatutin and Mikhailov~\cite{vm96}.

\section{The method of interaction graphs}\label{method}
One motivation for the present paper is to call attention to a new approach to proving central limit theorems for non-linear functions of independent random variables. Since most random variables can be so presented, the method has broad scope. The method is presented in \cite{cha08} in abstract form and used to solve a spatial statistics problem of Bickel; this involved `locally dependent' summands where the notion of local itself is determined from the data. A very different application is given in~\cite{chasound}. We hope that the applications in Sections \ref{proof} and \ref{future} will show the utility of this new approach.

The technique of defining `local neighborhoods' using the data was named `the method of interaction graphs' in \cite{cha08}. The method can be described very briefly as follows.  
Let $\xx$ be a set endowed with a sigma algebra and let $f:\xx^n \ra \rr$ be a measurable map, where $n\ge 1$ is a given positive integer. Suppose that $G$ is a map that associates to every $x\in \xx^n$ a simple graph $G(x)$ on $[n] := \{1,\ldots,n\}$. Such a map will be called a graphical rule on $\xx^n$. We will say that a graphical rule $G$ is symmetric if for any permutation $\pi$ of $[n]$ and any $(x_1,\ldots,x_n)\in \xx^n$, the set of edges in $G(x_{\pi(1)}, \ldots, x_{\pi(n)})$ is exactly
\[
\{\{\pi(i),\pi(j)\}: \{i,j\} \text{ is an edge of } G(x_1,\ldots,x_n)\}.
\]
For $m \ge n$, a symmetric graphical rule $G'$ on $\xx^m$ will be called an extension  of $G$ if for any $(x_1,\ldots,x_m) \in \xx^m$, $G(x_1,\ldots,x_n)$ is a subgraph of $G'(x_1,\ldots,x_m)$.

Now take any $x,x'\in \xx^n$. For each $i\in [n]$, let $x^i$ be the vector obtained by replacing $x_i$ with $x'_i$ in the vector $x$. For any two distinct elements $i$ and $j$ of~$[n]$, let $x^{ij}$ be the vector obtained by replacing $x_i$ with $x'_i$ and $x_j$ with~$x_j'$. We will say that the coordinates $i$ and $j$ are non-interacting for the triple $(f,x,x')$ if
\[
f(x)-f(x^j) = f(x^i) - f(x^{ij}).
\]
We will say that a graphical rule $G$ is an interaction rule for a function $f$ if for any choice of $x,x'$ and $i,j$, the event that $\{i,j\}$ is not an edge in the graphs $G(x)$, $G(x^i)$, $G(x^j)$, and $G(x^{ij})$ implies that $i$ and $j$ are non-interacting vertices for the triple $(f,x,x')$.  The following theorem implies central limit behavior if one can construct an interaction graph that has, with high probability, small maximum degree.
\begin{thm}[\cite{cha08}]\label{main}
Let $f:\xx^n \ra \rr$ be a measurable map that admits a symmetric interaction rule $G$. 
Let $X_1,X_2,\ldots $ be a sequence of i.i.d.\ $\xx$-valued random variables and let $X= (X_1,\ldots,X_n)$. Let $W := f(X)$ and $\sigma^2 := \var(W)$. Let $\xp = (\xp_1,\ldots,\xp_n)$ be an independent copy of $X$. For each $j$, define
\[
\Delta_j f(X) = W - f(X_1,\ldots,X_{j-1}, \xp_j, X_{j+1}, \ldots,X_n),
\]
and let $M = \max_j|\Delta_jf(X)|$. Let $G'$ be an extension of $G$ on $\xx^{n+4}$, and put
\[
\delta := 1+ \textup{degree of the vertex $1$ in } G'(X_1,\ldots,X_{n+4}).
\]
Then 
\[
\delta_W \le \frac{Cn^{1/2}}{\sigma^2}\ee(M^8)^{1/4} \ee(\delta^4)^{1/4} +  \frac{1}{2\sigma^3}\sum_{j=1}^n \ee|\Delta_j f(X)|^3,
\]
where $\delta_W$ is the Wasserstein distance between $(W-\ee(W))/\sigma$ and $N(0,1)$, and $C$ is a universal  constant.
\end{thm}

\section{Proof of Theorem \ref{clt}}\label{proof}
We now apply Theorem \ref{main} to prove Theorem \ref{clt}. Take $\mx = [0,1]^2$, and let $X_1,X_2,\ldots$ be independent uniformly distributed points on $\mx$. Let $X= (X_1,\ldots, X_n)$. Write each $X_i$ as a pair $(U_i,V_i)$. Define the $x$-rank of the point $X_i$ as the rank of $U_i$ among all the $U_j$'s, and the $y$-rank of the point $X_i$ as the rank of $V_i$ among all the $V_j$'s. More accurately, we should say ``$x$-rank of $X_i$ in $X$'' and ``$y$-rank of $X_i$ in $X$''.

Let $X_{(1)},\ldots,X_{(n)}$ be the $X_i$'s arranged according to their $x$-ranks, and let $X^{(1)},\ldots,X^{(n)}$ be the $X_i$'s arranged according to their $y$-ranks. Write $X_{(i)} = (U_{(i)}, V_{(i)})$ and $X^{(i)} = (U^{(i)}, V^{(i)})$. Let $\pi(i)$ be the $y$-rank of $X_{(i)}$. Then clearly $\pi$ is a uniform random permutation. Let $\sigma(i)$ be the $x$-rank of $X^{(i)}$. Then $X^{(i)} = X_{(\sigma(i))}$. Therefore 
\[
\pi(\sigma(i)) = \text{the $y$-rank of $X_{(\sigma(i))}$} = \text{the $y$-rank of $X^{(i)}$} = i\,.
\]
Thus, $\sigma = \pi^{-1}$. Let 
\begin{align*}
W = f(X) &:=  \sum_{i=1}^{n-1} 1_{\{\pi(i)> \pi(i+1)\}} + \sum_{i=1}^{n-1} 1_{\{\sigma(i)>\sigma(i+1)\}}\\
&= \sum_{i=1}^{n-1} 1_{\{V_{(i)} > V_{(i+1)}\}} + \sum_{i=1}^{n-1}1_{\{U^{(i)} > U^{(i+1)}\}}\,.
\end{align*}
Now suppose that one of the $X_i$'s is replaced by an independent copy $X_i'$. Then $W$ can change by at most $4$. Therefore in the notation of Theorem \ref{main}, $|\Delta_j f(X)|\le 4$ for every $j$ and hence $M\le 4$.  

For $x\in \mx^n$, define a simple graph $G(x)$ on $[n]$ as follows. For any $1\le i\ne j\le n$, let $\{i,j\}$ be an edge if and only if the $x$-rank of $x_i$ and the $x$-rank of $x_j$ differ by at most $1$, or the $y$-rank of $x_i$ and the $y$-rank of $x_j$ differ by at most $1$. This construction is clearly invariant under relabeling of indices, and hence $G$ is a symmetric graphical rule. 

Given $x,x'\in \mx^n$, let $x^i$, $x^j$ and $x^{ij}$ be defined as in the paragraphs preceding the statement of Theorem \ref{main}. Suppose that $\{i,j\}$ is not an edge in the graphs $G(x)$, $G(x^i)$, $G(x^j)$ and $G(x^{ij})$. Since the difference $f(x)-f(x^i)$ is determined only by those $x_k$'s such that $k$ is a neighbor of $i$ in  $G(x)$ or in $G(x^i)$, and the difference $f(x^j)-f(x^{ij})$ is determined only by those $x_k$'s such that $k$ is a neighbor of $i$ in $G(x^j)$ or in $G(x^{ij})$, the above condition implies that $i$ and $j$ are non-interacting for the triple $(f,x,x')$. Thus, $G$ is a symmetric interaction rule for $f$.

Next, define a graphical rule $G'$ on $\mx^{n+4}$ as follows. For any $x\in \mx^{n+4}$ and $1\le i\ne j\le n+4$, let $\{i,j\}$ be an edge if and only if the $x$-rank of $x_i$ and the $x$-rank of $x_j$ differ by at most $5$, or if the $y$-rank of $x_i$ and the $y$-rank of $x_j$ differ by at most $5$. Then $G'$ is a symmetric graphical rule. Since the addition of four extra points can alter the difference between the ranks of two points by at most $4$, any edge of $G(x_1,\ldots,x_n)$ is also an edge of $G'(x_1,\ldots,x_{n+4})$. Thus, $G'$ is an extension of $G$. 

Since the degree of $G'(x)$ is bounded by $10$ for any $x\in \mx^{n+4}$, the quantity $\delta$ of Theorem \ref{main} is bounded by $10$ in this example. The upper bounds on $\Delta_j f(X)$, $M$ and $\delta$ obtained above imply that for this $W$,
\[
\delta_W\le \frac{Cn^{1/2}}{\sigma^2} + \frac{Cn}{\sigma^3}\,,
\]
where $\sigma^2 = \var(W)$ and $\delta_W$ is the Wasserstein distance between $(W-\ee(W))/\sigma$ and $N(0,1)$. 

To complete the proof of Theorem \ref{clt}, we derive the given expression for $\var(T(\pi))$. We thank Walter Stromquist for the following elegant calculation. From classically available generating functions, it is straightforward to show that 
\[
\ee(D(\pi)) = \frac{n-1}{2}\,,\;\;\var(D(\pi)) = \frac{n+1}{12}\,,
\]
as given in Theorem \ref{descent} (see Fulman \cite{fulman04}). Consider, with obvious notation,
\[
\ee(D(\pi)D(\pi^{-1})) = \ee\biggl(\sum_{i=1}^{n-1} \sum_{j=1}^{n-1}D_i(\pi)D_j(\pi^{-1})\biggr)\,.
\]
For a given $\pi$, the $(n-1)^2$ terms in the sum may be broken into three types, depending on the size of the set $\{\pi(i), \pi(i+1)\} \cap \{j,j+1\}$, which may be $0$, $1$ or $2$. Carefully working out the expected value in each case, it follows that
\[
\ee(D(\pi)D(\pi^{-1})) =  \frac{(n-1)^2}{4} + \frac{n-1}{2n}\,.
\]
From this,
\[
\cov(D(\pi), D(\pi^{-1})) = \frac{n-1}{2n},\ \ \ \var(D(\pi)+D(\pi^{-1})) = 2\biggl(\frac{n+1}{12}\biggr) + 2\biggl(\frac{n-1}{2n}\biggr) = \frac{n+7}{6}-\frac{1}{n}\,.
\]
This completes the proof of Theorem \ref{clt}.
\vskip.1in
\noindent{\it \underline{Remarks.}} The argument used to prove Theorem \ref{clt} can be adapted to prove the joint limiting normality of $(D(\pi), D(\pi^{-1}))$. Fix real $(s,t)$ and consider $T_{s,t}(\pi)=s D(\pi)+t D(\pi^{-1})$. The same symmetric graph $G$ can be used. Because $s$, $t$ are fixed, the limiting normality follows as above. The limiting variance follows from $\cov(D(\pi), D(\pi^{-1}))$. We state the conclusion:
\begin{thm}\label{clt2}
For $\pi$ chosen uniformly from the symmetric group $S_n$, 
\[
\left(\frac{D(\pi)-\frac{n-1}{2}}{\sqrt{\frac{n+1}{12}}},\, \frac{D(\pi^{-1})-\frac{n-1}{2}}{\sqrt{\frac{n+1}{12}}}\right) \stackrel{d}{\longrightarrow} Z_2\,,
\]
where $Z_2$ is bivariate normal with mean $(0,0)$ and identity covariance matrix.
\end{thm}
Like Theorem \ref{clt}, Theorem \ref{clt2} can also be deduced as a corollary of the results in \cite{vatutin96}.

\section{Proof of Theorem \ref{genclt}}\label{proof2}
The proof goes exactly as the proof of Theorem \ref{clt}, except for slight differences in the definitions of the interaction graphs $G$ and $G'$. Recall the quantity $k$ from the statement of Theorem \ref{genclt}. With all the same notation as in the previous section, define the interaction graph $G$ as follows. For any $1\le i\ne j\le n$, let $\{i,j\}$ be an edge if and only if the $x$-rank of $x_i$ and the $x$-rank of $x_j$ differ by at most $k-1$, or the $y$-rank of $x_i$ and the $y$-rank of $x_j$ differ by at most $k-1$. It is easy to see that this is indeed an interaction graph for $F(\pi)+G(\pi^{-1})$. Define the extended graph $G'$ similarly, with $k-1$ replaced by $k+3$. Everything goes through as before, except that the quantity $M$ of Theorem \ref{main} is now bounded by $C_1 k$ (since the local components of $F$ and $G$ are bounded by $1$) and the quantity $\delta$ of Theorem \ref{main} is now bounded by $C_2 k$, where $C_1$ and $C_2$ are universal constants. The proof is now easily completed by plugging in these bounds in Theorem \ref{main}.

\section{Going further}\label{future}
There are two directions where the methods of this paper should allow natural limit theorems to be proved. The first is to develop a similar theory for {\it global} patterns in subsequences (for example triples $i<j<k$ with $\pi(i)<\pi(j)<\pi(k)$). Of course, inversions, the simplest case ($i<j$ with $\pi(i)>\pi(j)$), have a nice probabilistic and combinatorial theory, so perhaps the joint distribution of the number of occurrences of a fixed pattern in $\pi$ and $\pi^{-1}$ does as well. 

The second direction of natural generalization is to other finite reflection groups. There is a notion of `descents' --- the number of simple positive roots sent to negative roots. The question of number of number of descents in $x$ plus the number of descents in $x^{-1}$ makes sense and present methods should apply. In specific contexts, for example the hyperoctahedral group of signed permutations, these are concrete interesting questions. 

Lastly, it would be interesting to see if the dependence of the error bound on $k$ in Theorem \ref{genclt} can be improved, and to figure out what is the optimal dependence.

\end{document}